\documentclass[preprint,12pt]{elsarticle}

\usepackage{amsmath}
\usepackage{amssymb}

\newtheorem{lemma}{Lemma}
\newtheorem{proposition}{Proposition}
\newtheorem{definition}{Definition}
\newtheorem{theorem}{Theorem}
\newtheorem{corollary}{Corollary}

\def\Z{\mathbb{Z}}

\def\eps{\varepsilon}
\def\qed{$\blacksquare$}

\def\Oh{O}

\makeatletter
\def\pmod#1{\allowbreak\mkern10mu({\operator@font mod}\,\,#1)}

\makeatother

\journal{}

\begin{document}

\begin{frontmatter}

\title{On the second smallest prime non-residue}
\author{Kevin J. McGown\fnref{label}}
\fntext[label]{\emph{Current address:}  Department of Mathematics, Oregon State University, 368 Kidder Hall, Corvallis, OR 97331}
%% \author{Name\corref{cor1}\fnref{label2}}
\ead{mcgownk@math.oregonstate.edu}

\address{Department of Mathematics, University of California, San Diego,\\ 9500 Gilman Drive, La Jolla, CA 92093}

\begin{abstract}
Let $\chi$ be a non-principal Dirichlet character modulo a prime $p$.
Let $q_1<q_2$ denote the two smallest prime non-residues of $\chi$.
We give explicit upper bounds on $q_2$ that improve upon all known results.
We also provide a good upper estimate on the product $q_1 q_2$
which has an upcoming application to the study of norm-Euclidean Galois fields.
\end{abstract}

\begin{keyword}
%% keywords here, in the form: keyword \sep keyword
Dirichlet character \sep non-residues \sep power residues
%% MSC codes here, in the form: \MSC code \sep code
%% or \MSC[2008] code \sep code (2000 is the default)
\MSC[2010] Primary 11A15 \sep 11N25; Secondary 11A05
%\sep 11L26 \sep 11L40
\end{keyword}

\end{frontmatter}

\section{Introduction and Summary}\label{S:intro}
Let $\chi$ be a non-principal Dirichlet character modulo a prime $p$.
We call a positive integer $n$ a non-residue of $\chi$ if $\chi(n)\notin\{0,1\}$, and
denote by $q_1<q_2<\dots<q_n$ the $n$ smallest prime non-residues of $\chi$.
The question of putting an upper bound on $q_1$ is a classical problem
which goes all the way back to the study of the least quadratic non-residue.

The literature on this problem is extensive and we will not review it here except to say that 
the work of Burgess in the 1960's significantly advanced existing knowledge on this matter.
Burgess' famous character sum estimate (see~\cite{burgess:1962})
implies that $q_n=\Oh(p^{1/4+\eps})$ for all $n$.\footnote
{
The $\Oh$ constant here depends upon $\eps$ and $n$; see~\cite{hudson:1983} for more detail.
}
For the case of $q_1$, one can apply the ``Vinogradov trick''
(see~\cite{vinogradov:1918, vinogradov:1927, vinogradov:1927b})
to Burgess' result, which
gives the stronger bound of
$q_1=\Oh(p^{\frac{1}{4\sqrt{e}}+\eps})$
(see~\cite{burgess:1962}).
%Additionally, various authors have given improvements to this $\Oh$-bound on $q_1$ which depend on the order of~$\chi$.

Making these results explicit with constants of a reasonable magnitude
turns out to be difficult, and often times it is results of this nature that
one requires in application.
In this paper, we will restrict ourselves to the study of $q_1$ and $q_2$,
and we will only be interested in bounds which are completely explicit and independent of the order of
$\chi$.\footnote{In Corollary~\ref{C:3} we do assume that $\chi$ has odd order, but we emphasize that none of our constants depend upon the order of $\chi$.}

The best known explicit bound on $q_1$ was given by Norton (see~\cite{norton:1971})
by applying Burgess' method (see~\cite{burgess:1962, burgess:1963}) with some modifications.
%Namely, he shows that $q_1<4.7\,p^{1/4}\log p$,
%and moreover, that the constant can be improved to $3.9$ when the order of $\chi$ and $(p-1)/2$ have a common factor.
\begin{theorem}[Norton]\label{T:norton}
  Suppose that $\chi$ is a non-principal Dirichlet character modulo a prime~$p$, and that
  $q_1$ is the smallest (prime) non-residue of $\chi$. 
  Then $q_1<4.7\,p^{1/4}\log p$, and moreover, 
  the constant can be improved to $3.9$ when the order of $\chi$ and $(p-1)/2$ have a common factor.
\end{theorem}
We prove the following theorem, which can be viewed as a generalization of Norton's result
but with a slightly larger constant.
\begin{theorem}\label{T:1}
  Fix a real constant $p_0\geq 10^7$.
  There exists an explicit constant $C$ (see Table~\ref{Table:1})
  such that if
  $\chi$ is a non-principal Dirichlet character modulo a prime $p\geq p_0$
  and $u$ is a prime with $u\geq e^2\log p$,
  then there exists $n\in\Z^+$ with $(n,u)=1$, $\chi(n)\neq 1$, and
  $$
    n<C\,p^{1/4}\log p
    \,.
  $$  
  %\vspace{-2ex}
\begin{table}[t]
\centering
\begin{tabular}{| l | r |}
\hline
$p_0$ & $C\phantom{123}$ \\
\hline
$10^{7}$ & $11.0421$\\
$10^{8}$ & $8.2760$\\
$10^{9}$ & $7.2906$\\
$10^{10}$ & $6.8121$\\
$10^{11}$ & $6.5496$\\
$10^{12}$ & $6.3964$\\
$10^{13}$ & $6.3033$\\
\hline
\end{tabular}
\qquad
\begin{tabular}{| l | r |}
\hline
$p_0$ & $C\phantom{123}$ \\
\hline
$10^{14}$ & $6.2452$\\
$10^{15}$ & $6.2077$\\
$10^{16}$ & $6.1829$\\
$10^{17}$ & $6.1659$\\
$10^{18}$ & $6.1536$\\
$10^{19}$ & $6.1445$\\
$10^{20}$ & $6.1374$\\
\hline
\end{tabular}
\caption{Values of $C$ for various choices of $p_0$\label{Table:1}}
\end{table}
\end{theorem}
Provided that $q_1$ is not too small, the above theorem
immediately gives an explicit bound on $q_2$.
%The above theorem immediately implies the following result
%which gives an explicit bound on $q_2$ which is almost as good as
%Norton's bound on $q_1$, provided that $q_1$ is ``not too small''.
\begin{corollary}\label{C:1}
Fix a real constant $p_0\geq 10^7$.
Let $\chi$ be a non-principal Dirichlet character modulo a prime $p\geq p_0$.
%Suppose $q_1<q_2$ are the two smallest primes such that $\chi(q_1),\chi(q_2)\neq 1$.
Suppose that $q_1<q_2$ are the two smallest prime non-residues of $\chi$. 
If $q_1>e^2\log p$, then
$$
  q_2<C\,p^{1/4}\log p
  \,,
$$
where the constant $C$ is the same constant as in the statement
of Theorem~\ref{T:1} (see Table~\ref{Table:1}).
\end{corollary}

Using a lemma of Hudson and an explicit result of the author on consecutive non-residues,
we can remove the restriction on $q_1$
%at the price of increasing the size of the constant.
for a small price.
\begin{corollary}\label{C:2}
Let $\chi$ be a non-principal Dirichlet character modulo a prime $p\geq 10^{19}$.
%Suppose $q_1<q_2$ are the two smallest primes such that $\chi(q_1),\chi(q_2)\neq 1$.
Suppose that $q_1<q_2$ are the two smallest prime non-residues of $\chi$.
Then
$$
  q_2<53\,p^{1/4}(\log p)^2
  \,.
$$
\end{corollary}
The value $q_2$ has not been as extensively studied as $q_1$, and it appears that
prior to now, the best explicit bound was essentially
$q_2\leq c\,p^{2/5}$ for some absolute constant $c$ (see~\cite{brauer:1931, brauer:1940, whyburn, hudson:1973}).
Corollary~\ref{C:2} constitutes an explicit bound on $q_2$ which even improves slightly on
the best known $\Oh$-bound of $p^{1/4+\eps}$.

For the application the author has in mind to norm-Euclidean Galois fields
(see~\cite{mcgown:euclidean}),
the following corollary is more useful.
\begin{corollary}\label{C:3}
Let $\chi$ be a non-principal Dirichlet character modulo a prime $p\geq 10^{18}$
having odd order.
%Suppose $q_1<q_2$ are the two smallest primes such that $\chi(q_1),\chi(q_2)\neq 1$.
Suppose that $q_1<q_2$ are the two smallest prime non-residues of $\chi$.
Then
$$
  q_1 q_2<24\,p^{1/2}(\log p)^2
  \,.
$$
\end{corollary}

%As before, the constant $C'$ only depends upon $p_0$.
%and the following is a table with admissible values of $C'$ for various choices of $p_0$:

%\begin{table}[H]
%\centering
%\caption{Values of $C'$ for various choices of $p_0$}
%\begin{tabular}{c c}
%\begin{tabular}{| l | r |}
%\hline
%$p_0$ & $C'\phantom{123}$ \\
%\hline
%$10^{7}$ & $43.0892$\\
%$10^{8}$ & $32.2831$\\
%$10^{9}$ & $28.4365$\\
%$10^{10}$ & $26.5688$\\
%$10^{11}$ & $25.5443$\\
%$10^{12}$ & $24.9464$\\
%$10^{13}$ & $24.5833$\\
%\hline
%\end{tabular}
%\qquad
%\begin{tabular}{| l | r |}
%\hline
%$p_0$ & $C'\phantom{123}$ \\
%\hline
%$10^{14}$ & $24.3563$\\
%$10^{15}$ & $24.2105$\\
%$10^{16}$ & $24.1134$\\
%$10^{17}$ & $24.0471$\\
%$10^{18}$ & $23.9995$\\
%$10^{19}$ & $23.9636$\\
%$10^{20}$ & $23.9359$\\
%\hline
%\end{tabular}
%\end{tabular}
%\end{table}
%

%%%%%%%%%%%%%%%%%%%%%%%%%%%%%%%%%%%%%%%%%%%%%%%%%%%%%%%%%%%%%%%%%%%%%%%%%%%%%%%%
%%%%%%%%%%%%%%%%%%%%%%%%%%%%%%%%%%%%%%%%%%%%%%%%%%%%%%%%%%%%%%%%%%%%%%%%%%%%%%%%
%%%%%%%%%%%%%%%%%%%%%%%%%%%%%%%%%%%%%%%%%%%%%%%%%%%%%%%%%%%%%%%%%%%%%%%%%%%%%%%%

\section{Outline of the Proof}\label{S:outline}
We will establish our results using a generalization of Burgess' method.
The approach will be similar to a previous paper of the author (see~\cite{mcgown:consecutive}),
but it will be sufficiently different as these results do not follow from the aforementioned
ones or vice versa.
The main idea behind Burgess' method is to combine upper and lower bounds for the following sum:
\begin{definition}
%If $\chi$ is a Dirichlet character modulo $p$ and $h,r\in\Z^+$, then we define
If $h,r\in\Z^+$ and $\chi$ is a Dirichlet character modulo $p$, then we define
$$
 S(\chi,h,r):= \sum_{x=0}^{p-1}
  \left|
  \sum_{m=1}^h
  \chi(x+m)
  \right|^{2r}
  \,.
$$
\end{definition}

We will use the following lemma, proven in~\cite{mcgown:consecutive},
which is a slight improvement on
Lemma 2 of~\cite{burgess:1962}.
\begin{lemma}\label{L:1C}
Suppose $\chi$ is any non-principal Dirichlet character to the prime modulus $p$.
If $r,h\in\Z^+$, then
$$
 S(\chi,h,r)
  <
  \frac{1}{4}(4r)^rp h^r
  +
  (2r-1)p^{1/2}h^{2r}
  \,.
$$  
\end{lemma}

Apart from the use of Lemma~\ref{L:1C}, the proofs of Theorem~\ref{T:1} and Corollary~\ref{C:1} are
completely self-contained; in particular, they do not rely on Theorem~\ref{T:norton}.
However, the derivation of Corollary~\ref{C:2} will use Theorem 1.2 of~\cite{mcgown:consecutive},
and the derivation of Corollary~\ref{C:3} will use Theorem~\ref{T:norton}
and an explicit version of the P\'olya--Vinogradov inequality given in~\cite{bachman.rachakonda}.

The meat of the proof of our results is to give
a lower bound on $S(\chi,h,r)$, under some extra conditions on the involved parameters.
In \S\ref{S:prop} we prove the following:
\begin{proposition}\label{P:q2UB}
  Let $h,r,u\in\Z^+$ with $u$ prime and $h\leq u$.
  Suppose that $\chi$ is a Dirichlet character modulo a prime $p\geq 5$
  such that $\chi(n)=1$ for all $n\in[1, H]$ satisfying $(n,u)=1$.
  Assume $2h<H\leq(2hp)^{1/2}$ and set $X:=H/(2h)>1$.
  Then
$$
  S(\chi,h,r)\geq 
   \frac{6}{\pi^2}
  (1-u^{-1})h(h-2)^{2r}X^2 f(X,u)
  \,.
$$
For each fixed $u$ we have $f(X,u)\to 1$ as $X\to\infty$;
the function $f(X,u)$ is explicitly defined in Lemma~\ref{L:phisum}.
\end{proposition}

Combining Lemma~\ref{L:1C} and Proposition~\ref{P:q2UB} with a careful choice of the parameters
$h$ and $r$
gives our main result from which Theorem~\ref{T:1} follows:

\begin{theorem}\label{T:2}
  Suppose that $\chi$ is a non-principal Dirichlet character modulo a prime $p\geq 10^7$,
  and that $u$ is a prime with $u\geq e^2\log p$.
  Suppose $\chi(n)=1$ for all $n\in[1,H]$ with $(n,u)=1$.
  If
  $$
    H\leq (2e^2\log p-2)^{1/2}p^{1/2}
  \,,
  $$
  then
  $$
    H\leq
    Kg(p)\,p^{1/4}\log p
    \,,
  $$
  where
  $$
  K=\frac{\pi e}{\sqrt{2}}
  \approx
  %6.038504
  6.0385
  $$
  and
  $$
    g(p)=
    \sqrt{
    \frac{\left(1+\frac{4}{3\log p}\right)}
    {
      \left(1-\frac{1}{e^2\log p}\right)f\left(\frac{Kp^{1/4}}{2e^2},89\right)
    }
    }
    \,.
  $$
  The function $g(p)$ is positive and decreasing for $p\geq 10^7$, with $g(p)\to 1$ as $p\to\infty$.
  The function $f(X,u)$ is defined in Lemma~\ref{L:phisum}.
\end{theorem}

The proofs of Theorems~\ref{T:1} and~\ref{T:2} are carried out in \S\ref{S:theorem}.
Finally in \S\ref{S:corollaries} we derive Corollaries~\ref{C:1}, \ref{C:2}, and \ref{C:3}.

%%%%%%%%%%%%%%%%%%%%%%%%%%%%%%%%%%%%%%%%%%%%%%%%%%%%%%%%%%%%%%%%%%%%%%%%%%%%%%%%
%%%%%%%%%%%%%%%%%%%%%%%%%%%%%%%%%%%%%%%%%%%%%%%%%%%%%%%%%%%%%%%%%%%%%%%%%%%%%%%%
%%%%%%%%%%%%%%%%%%%%%%%%%%%%%%%%%%%%%%%%%%%%%%%%%%%%%%%%%%%%%%%%%%%%%%%%%%%%%%%%

\section{Proof of Proposition~\ref{P:q2UB}}\label{S:prop}

%The meat of the proof of Theorem~\ref{T:q2main} is to give
%a lower bound on $S(\chi,h,r)$, under some extra conditions on the involved parameters.
The idea is to locate a large number
of disjoint intervals on which $\chi$ is ``almost constant.''
For the remainder of this section $p$ will denote a prime with $p\geq 5$, and
$h,H$ will denote positive integers.
The following are the intervals that will be of interest to us:

\begin{definition}
  For integers with $0\leq t<q$, 
  we define the intervals
  \begin{eqnarray*}
    &&
    \mathcal{I}(q,t)=\left(\frac{pt}{q},\frac{H+pt}{q}\right]
    \,,
    \quad
        \mathcal{I}(q,t)^\star=\left(\frac{pt}{q},\frac{H+pt}{q}-h\right]
        \,,
    \\
    &&
    \mathcal{J}(q,t)=\left[-\frac{H+pt}{q},-\frac{pt}{q}\right)
    \,,
    \quad
    \mathcal{J}(q,t)^\star=\left[-\frac{H+pt}{q},-\frac{pt}{q}-h\right)
    \,.
  \end{eqnarray*}
\end{definition}

We note that the intervals $\mathcal{I}(q,t)^\star$, $\mathcal{J}(q,t)^\star$ might be empty.
In fact, they are non-empty exactly when $h<H/q$, which will always be the case
whenever we employ them.

%%%%%%%%%%%%%%%%%%%%%%%%%%%%%%%%%%%%%%%%%%%%%%%%%%%%%%%%%%%%%%%%%%%%%%%%

\begin{lemma}\label{L:IJdisjoint}
  Let $X>1$ be a real number and suppose $XH<p$.  Then the intervals $\mathcal{I}(q,t)$ where $0\leq t<q\leq X$ with $(t,q)=1$ are disjoint, and
  similarly for $\mathcal{J}(q,t)$.
\end{lemma}

\noindent\textbf{Proof.}
If $\mathcal{I}(q_1,t_1)$ and $\mathcal{I}(q_2,t_2)$ intersect, then we have:
%\begin{eqnarray*}
%&&
%(t_1 p-H)q_1\leq (t_2p+H)/q_2
%\\
%&&
%(t_2 p-H)q_2\leq (t_1p+H)/q_1
%\end{eqnarray*}
\begin{eqnarray*}
&&
pt_1/q_1\leq(H+pt_2)/q_2
\\
&&
pt_2/q_2\leq(H+pt_1)/q_1
\end{eqnarray*}
It follows that
$$
  |t_1q_2-t_2 q_1|
  \leq
  \frac{XH}{p}
  <1
  \,;
$$
whence $t_1 q_2=t_2 q_1$ which implies
$t_1=t_2$, $q_1=q_2$.
(When $t_1=t_2=0$, the condition $(q_1,t_1)=(q_2=t_2)=1$ forces
$q_1=q_2=1$, so the argument goes through in this case as well.)
The proof for the intervals $\mathcal{J}(q,t)$ is the same.~\qed

%%%%%%%%%%%%%%%%%%%%%%%%%%%%%%%%%%%%%%%%%%%%%%%%%%%%%%%%%%%%%%%%%%%%%%%%

\begin{lemma}\label{L:chisumhminus2}
  Let $h,u\in\Z^+$ with $u$ prime and $h\leq u$.
  Suppose that $\chi$ is a Dirichlet character modulo $p$ such that $\chi(n)=1$ for all $n\in[1,H]$ with $(n,u)=1$.
  If $z\in\mathcal{I}(q,t)^\star\cup\mathcal{J}(q,t)^\star$ and $(q,u)=1$, then
  %then the values $\chi(z+n)$ for $n=0,\dots, h-1$ are all equal except for possibly one value of $n$, and
  %likewise for $\mathcal{J}(q,t)^\star$.
  $$
  \left|
 % \sum_{m=1}^h
  \sum_{m=0}^{h-1}
  \chi(z+m)
  \right|
  \geq
  h-2
  \,.
  $$
\end{lemma}

\noindent\textbf{Proof.}
We note that by hypothesis $\mathcal{I}(q,t)^\star\cup\mathcal{J}(q,t)^\star\neq\emptyset$
and hence $h<H/q$.
First suppose $z\in\mathcal{I}(q,t)^\star$.
We will show that
the values $\chi(z+n)$ for
%$n=1,\dots, h$
$n=0,\dots,h-1$
are all equal except for possibly one value of $n$.
%, and likewise for $\mathcal{J}(q,t)^\star$.
This will immediately give the result upon application of the
triangle inequality.

For
%$n=1,\dots,h$,
$n=0,\dots,h-1$,
we have
%$z\in\mathcal{I}(q,t)^\star$, which implies
$z+n\in\mathcal{I}(q,t)$ and hence $q(z+n)-pt\in(0,H]$.
Provided $u$ does not divide $q(z+n)-pt$, we have
$$
  \chi(z+n)=\overline{\chi}(q)\chi(q(z+n))=\overline{\chi}(q)\chi(q(z+n)-pt)=\overline{\chi}(q)
  \,.
$$
But if $u$ divides $q(z+n)-pt$ for two distinct values of $n$, say
$n_1$ and $n_2$, we find that $u$ divides $q(n_1-n_2)$.  Since $(u,q)=1$, we conclude
that $u$ divides $n_1-n_2$ and hence $|n_1-n_2|\geq u$.  This leads to
$h\leq u\leq |n_1-n_2|\leq h-1$, a contradiction.
The proof for $z\in\mathcal{J}(q,t)^\star$ is similar.
\qed

%%%%%%%%%%%%%%%%%%%%%%%%%%%%%%%%%%%%%%%%%%%%%%%%%%%%%%%%%%%%%%%%%%%%%%%%

\begin{lemma}\label{L:SR1}
  Suppose that $X>1$ is a real number and $u\in\Z^+$ is prime.
  Then
  $$
    %\sum_{\footnotesize{\begin{array}{cc}n\leq X\\[-1ex](n,u)=1\end{array}}}
    \sum_{\substack{n\leq X\\(n,u)=1}}
    n
    =\frac{\left(1-u^{-1}\right)}{2}X^2+\theta_{X,u}X
    \,,
  $$
  where the sum is taken over positive integers and $\theta_{X,u}$ denotes a real number,
  depending on $X$ and $u$, that belongs to
  the interval $(-1,1)$.
\end{lemma}

\noindent\textbf{Proof.}
For any $Y>0$ we have
$$
  \sum_{n\leq Y}
  n
  =\frac{\lfloor Y\rfloor(\lfloor Y\rfloor+1)}{2}
  \,.
$$
Upon an application of the obvious inequality $Y-1<\lfloor Y\rfloor\leq Y$, we obtain
the identity
$$
  \sum_{n\leq Y}
  n
  =\frac{Y^2}{2}
  +
  \frac{Y}{2}
  \theta_Y
  \,,
$$
where $\theta_Y\in(-1,1]$.
Now we write
\begin{eqnarray*}
    %\sum_{\footnotesize{\begin{array}{cc}n\leq X\\[-1ex](n,u)=1\end{array}}}
    \sum_{\substack{n\leq X\\(n,u)=1}}
    n
    &=&
    \sum_{n\leq X}n
    -
    u\sum_{n\leq X/u}
    n
    \\
    &=&
    \frac{X^2}{2}(1-u^{-1})+\frac{X}{2}(\theta_X-\theta_{X/u})
    \,,
\end{eqnarray*}
and observe that
$$
  -2<\theta_X-\theta_{X/u}<2
  \,.
$$
The result follows.
\qed

%%%%%%%%%%%%%%%%%%%%%%%%%%%%%%%%%%%%%%%%%%%%%%%%%%%%%%%%%%%%%%%%%%%%%%%%

\begin{lemma}\label{L:phisum}
Suppose $X>1$ and $u\in\Z^+$ is prime.  Then
$$
%  \sum_{\footnotesize{\begin{array}{cc}1\leq q\leq X\\(q,u)=1\end{array}}}
%  \sum_{\footnotesize{\begin{array}{cc}0\leq t<q\\(t,q)=1\end{array}}}
%  1
  %\#\{(q,t)\in\Z^2\mid 0\leq t<q\leq X,\,(q,tu)=1\}
  %\sum_{\footnotesize{\begin{array}{cc}1\leq q\leq X\\[-1ex](q,u)=1\end{array}}}
  \sum_{\substack{1\leq q\leq X\\(q,u)=1}}
  \phi(q)
  \;\geq\;
  \frac{3}{\pi^2}
  (1-u^{-1})X^2f(X,u)
  \,,
$$
where
$$
  f(X,u)=1-\frac{\pi^2}{3}
  \left(
  \frac{1}{2X^2}+\frac{1}{2X}+\frac{1}{1-u^{-1}}\cdot\frac{1+\log X}{X}
  \right)
  \,.
$$
\end{lemma}

\noindent\textbf{Proof.}
First we observe:
\begin{eqnarray*}
%\sum_{\footnotesize{\begin{array}{cc}1\leq q\leq X\\[-1ex](q,u)=1\end{array}}}
\sum_{\substack{1\leq q\leq X\\(q,u)=1}}
\phi(q)
&=&
%\sum_{\footnotesize{\begin{array}{cc}1\leq q\leq X\\[-1ex](q,u)=1\end{array}}}
\sum_{\substack{1\leq q\leq X\\(q,u)=1}}
\sum_{m|q}
\frac{q}{m}\;\mu(m)
\\[1ex]
&=&
\sum_{\substack{1\leq m\leq X\\(m,u)=1}}
\mu(m)
\sum_{\substack{1\leq r\leq X/m\\(r,u)=1}}
r
\end{eqnarray*}
Applying Lemma~\ref{L:SR1} to the above gives:
\begin{eqnarray*}
  &&
  %\sum_{\footnotesize{\begin{array}{cc}1\leq q\leq X\\[-1ex](q,u)=1\end{array}}}
  \sum_{\substack{1\leq q\leq X\\(q,u)=1}}
  \phi(q)
  \;\;=\;\;
  \\
  && 
  \quad
  \frac{X^2}{2}\left(1-u^{-1}\right)
  \left(
  \sum_{\substack{1\leq m\leq X\\(m,u)=1}}
  \frac{\mu(m)}{m^2}
    \right)
  %\quad
  +\;
  X
  \left(
  %\sum_{\footnotesize{\begin{array}{cc}1\leq m\leq X\\[-1ex](m,u)=1\end{array}}}
  \sum_{\substack{1\leq m\leq X\\(m,u)=1}}
  \frac{\mu(m)}{m}
  \,
  \theta_{X/m,u}  
      \right)
\end{eqnarray*}
Now we use the bounds: % (see the proof of Lemma~3.1 in~\cite{mcgown:euclidean}):
$$
    %\sum_{\footnotesize{\begin{array}{cc}1\leq m\leq X\\[-1ex](m,u)=1\end{array}}}
    \sum_{\substack{1\leq m\leq X\\(m,u)=1}}
    \frac{\mu(m)}{m^2}
    \geq
    \frac{6}{\pi^2}
    -
    \frac{1}{X^2}
    -
    \frac{1}{X}
    \,,
$$
$$
%\left|
  %\sum_{\footnotesize{\begin{array}{cc}1\leq m\leq X\\[-1ex](m,u)=1\end{array}}}
  \vrule \
    \sum_{\substack{1\leq m\leq X\\(m,u)=1}}^{\phantom{n}}
  \frac{\mu(m)}{m}
  \,
  \theta_{X/m,u}  
  \ \vrule
  \leq
  %\vrule \
  %\sum_{\substack{1\leq m\leq X\\(m,u)=1}}^{\phantom{n}}
  \sum_{1\leq m\leq X}
  \frac{1}{m}
  %\ \vrule
 % \right|
  \leq
  1+\log X
$$
The result follows from an application of the triangle inequality and some rearrangement.
\qed

%%%%%%%%%%%%%%%%%%%%%%%%%%%%%%%%%%%%%%%%%%%%%%%%%%%%%%%%%%%%%%%%%%%%%%%%
%\vspace{1ex}
%Finally, we are ready to give the lower bound we have alluded to.
%
%\begin{proposition}\label{P:q2UB}
%  Let $h,r,u\in\Z^+$ with $u$ prime and $h\leq u$.
%  Suppose that $\chi$ is a Dirichlet character modulo $p$ such that $\chi(n)=1$ for all $n\in[1, H]$ satisfying $(n,u)=1$.
%  Assume $2h<H\leq(2hp)^{1/2}$ and set $X:=H/(2h)>1$.
%  Then
%$$
%  S(\chi,h,r)\geq 
%   \frac{6}{\pi^2}
%  (1-u^{-1})h(h-2)^{2r}X^2 f(X,u)
%  \,.
%$$
%The function $f(X,u)$ is defined in Lemma~\ref{L:phisum}.
%\end{proposition}

\vspace{1ex}
\noindent\textbf{Proof of Proposition~\ref{P:q2UB}.}
We begin by noting that $H/q\geq H/X=2h$.
Using Lemma~\ref{L:IJdisjoint} and Lemma~\ref{L:chisumhminus2}  we have:
\begin{eqnarray*}
  S(\chi,h,r)
  &=&
  \sum_{x=0}^{p-1}\left|\sum_{m=0}^{h-1}\chi(x+m)\right|^{2r}
  \\[1ex]
  &\geq&
  \sum_{\substack{0\leq t<q\leq X\\(q,u)=(q,t)=1}}
  \sum_{z\in\mathcal{I}_{q,t}^\star\cup\mathcal{J}_{q,t}^\star}
  \left|
  \sum_{m=0}^{h-1}
  \chi(z+m)
  \right|^{2r}
  \\[1ex]
  &\geq&
  \sum_{\substack{0\leq t<q\leq X\\(q,tu)=1}}
  2\left(\frac{H}{q}-h\right)(h-2)^{2r}
  \\[1ex]
  &\geq&
  \sum_{\substack{0\leq t<q\leq X\\(q,tu)=1}}
  2h(h-2)^{2r}
 \\[1ex]
  &=&
  2h(h-2)^{2r}
  \sum_{\substack{1\leq q\leq X\\(q,u)=1}}
  \phi(q)
  %\,.
\end{eqnarray*}
%\begin{eqnarray*}
%  S(\chi,h,r)
%  &=&
%  \sum_{x=0}^{p-1}\left|\sum_{m=0}^{h-1}\chi(x+m)\right|^{2r}
%  \\[1ex]
%  &\geq&
%  \sum_{\footnotesize{\begin{array}{cc}0\leq t<q\leq X\\[-1ex](q,u)=(q,t)=1\end{array}}}
%    \sum_{\footnotesize{\begin{array}{c}z\in\mathcal{I}_{q,t}^\star\cup\mathcal{J}_{q,t}^\star\end{array}}}
%  \left|
%  \sum_{m=0}^{h-1}
%  \chi(z+m)
%  \right|^{2r}
%  \\[1ex]
%  &\geq&
%  \sum_{\footnotesize{\begin{array}{cc}0\leq t<q\leq X\\[-1ex](q,tu)=1\end{array}}}
%  2\left(\frac{H}{q}-h\right)(h-2)^{2r}
%  \\[1ex]
%  &\geq&
%  \sum_{\footnotesize{\begin{array}{cc}0\leq t<q\leq X\\[-1ex](q,tu)=1\end{array}}}
%  2h(h-2)^{2r}
%  \,,
%  \end{eqnarray*}
%and hence
%\begin{eqnarray*}
%  S(\chi,h,r)
%  &\geq&
%    \sum_{\footnotesize{\begin{array}{cc}0\leq t<q\leq X\\[-1ex](q,tu)=1\end{array}}}
%  2h(h-2)^{2r}
% \\[1ex]
%  &=&
%  2h(h-2)^{2r}
%  \sum_{\footnotesize{\begin{array}{cc}1\leq q\leq X\\[-1ex](q,u)=1\end{array}}}
%  \phi(q)
%  \,.
%\end{eqnarray*}
Now the result follows from Lemma~\ref{L:phisum}.
\qed

%%%%%%%%%%%%%%%%%%%%%%%%%%%%%%%%%%%%%%%%%%%%%%%%%%%%%%%%%%%%%%%%%%%%%%%%

\section{Proofs of the Theorems}\label{S:theorem}

Before launching the proof of Theorem~\ref{T:2}, we establish the following simple convexity result:

\begin{lemma}\label{L:convex}
Suppose $h,r\geq 1$.  We have the following implications:
%$$
%  \frac{1}{2h}\left(\frac{4r}{h-2}\right)^r
%  \leq
%  \frac{1}{h+1}
%  \left(\frac{4r}{h+1}\right)^r
%  \,.
%$$
\begin{eqnarray*}
 h\geq 6r+5
  &\Longrightarrow&
  \frac{1}{2h}\left(\frac{4r}{h-2}\right)^r
  \leq
  \frac{1}{h+1}
  \left(\frac{4r}{h+1}\right)^r
  \\[1ex]
  h\geq 16r+2
  &\Longrightarrow&
  \left(\frac{h}{h-2}\right)^r<\frac{7}{6}
  \\[1ex]
  h\geq 2r-1
  &\Longrightarrow&
  \frac{2r-1}{h}
  \leq
  \frac{2r}{h+1}
\end{eqnarray*}
\end{lemma}

\noindent\textbf{Proof.}
By the convexity of the logarithm,
we have $\log t\geq (2\log 2)(t-1)$ for all $t\in[1/2,1]$.
Applying this, together with the hypothesis that $6(r+1)\leq h+1$, we get
$$
  \log\left(\frac{h-2}{h+1}\right)
  \geq
  -\frac{6\log 2}{h+1}
  \geq
  -\frac{\log 2}{r+1}
  \,.
$$
This yields
$$
  \frac{1}{2}
  \leq
  \left(\frac{h-2}{h+1}\right)^{r+1}
%  \leq
%  \frac{h}{h+1}
%  \left(
%  \frac{h-2}{h+1}
%  \right)^{r+1}
  \,,
$$
and first implication follows.
For the proof of the second implication, we observe (again by convexity) that
$\log t\leq t-1$ for all $t$ and hence
$$
  r\log\left(\frac{h}{h-2}\right)\leq\frac{2r}{h-2}\leq\frac{1}{8}
  \,;
$$
this leads to
$$
  \left(\frac{h}{h-2}\right)^r\leq\exp\left(\frac{1}{8}\right)<\frac{7}{6}
  \,.
$$
The third implication is trivial.
\qed

%\vspace{1ex}

%The following is the main result of the paper, from which Theorem~\ref{T:q2main} will follow:
%
%%%%%%%%%%%%%%%%%%%%%%%%%%%%%%%%%%%%%%%%%%%%%%%%%%%%%%%%%%%%%%%%%%%%%%%%%
%
%\begin{theorem}\label{T:q2maintech}
%  Suppose that $\chi$ is a non-principal Dirichlet character modulo $p\geq 10^7$,
%  and that $u$ is a prime with $u\geq e^2\log p$.
%  Suppose $\chi(n)=1$ for all $n\in[1,H]$ with $(n,u)=1$.
%  If
%  $$
%    H\leq (2e^2\log p-2)^{1/2}p^{1/2}
%  \,,
%  $$
%  then
%  $$
%    H\leq
%    Kg(p)\,p^{1/4}\log p
%    \,,
%  $$
%  where
%  $$
%  K=\frac{\pi e}{\sqrt{2}}
%  \approx
%  %6.038504
%  6.0385
%  $$
%  and
%  $$
%    g(p)=
%    \sqrt{
%    \frac{\left(1+\frac{4}{3\log p}\right)}
%    {
%      \left(1-\frac{1}{e^2\log p}\right)f\left(\frac{Kp^{1/4}}{2e^2},89\right)
%    }
%    }
%    \,.
%  $$
%  The function $g(p)$ is positive and decreasing for $p\geq 10^7$, with $g(p)\to 1$ as $p\to\infty$.
%  The function $f(X,u)$ is defined in Lemma~\ref{L:phisum}.
%\end{theorem}
%

\vspace{2ex}
\noindent\textbf{Proof of Theorem~\ref{T:2}.}
First, we may assume $H\geq Kp^{1/4}\log p$, or else there is nothing to prove.
We set $h=\lfloor A\log p\rfloor$, $r=\lfloor B\log p\rfloor$ with $A=e^2$, $B=1/4$
and verify that $r,h$ satisfy all three conditions in Lemma~\ref{L:convex}.
The constants $A$ and $B$ were chosen to minimize the quantity $AB$ subject to the constraint
$A\geq 4B\exp(1/(2B))$.

One verifies that $Kp^{1/4}> 28e^2$ for $p\geq 10^7$ and hence $H> 28h$.
%We set $X:=H/(2h)$ and observe that $X>7$.
We set $X:=H/(2h)$ and observe that we have the a priori lower bound
$$
  X=
  \frac{H}{2h}
  \geq
  \frac{Kp^{1/4}\log p}{2e^2\log p}
  =
  \frac{Kp^{1/4}}{2e^2}
  \,,
$$
and, in particular, $X>14$ from the previous sentence.
Since $p\geq 10^5$ and $e^2\log (10^5)\approx 85.1$, we know $u\geq 89$ and hence $f(X,u)\geq f(X,89)$.
For notational convenience, we will write $f(X):=f(X,89)$.

Combining Lemma \ref{L:1C} and Proposition \ref{P:q2UB}, we obtain
$$
  \frac{6}{\pi^2}\left(1-u^{-1}\right)h(h-2)^{2r}\left(\frac{H}{2h}\right)^2f(X)
  \;\leq\;
  \frac{1}{4}(4r)^r ph^r + (2r-1)p^{1/2}h^{2r}
  \,.
$$
Rearranging the above and applying Lemma~\ref{L:convex} gives
\begin{eqnarray}
&&
\nonumber
\frac{6}{\pi^2}\left(1-u^{-1}\right)H^2f(X)
\\
\nonumber
&&\qquad
\leq
4h^2p^{1/2}
\left[
\frac{1}{4h}\left(\frac{4r}{h-2}\right)^r\left(\frac{h}{h-2}\right)^r p^{1/2}+\frac{2r-1}{h}\left(\frac{h}{h-2}\right)^{2r}\right]
\\[1ex]
\label{E:explain1}
&&\qquad
\leq
4h^2p^{1/2}
\left[
\frac{1}{h+1}\left(\frac{4r}{h+1}\right)^r p^{1/2}+\frac{3r}{h+1}
\right]
\,.
\end{eqnarray}
Plugging in our choices of $r,h$ and using the fact that
$$
  A\geq 4B\exp\left(\frac{1}{2B}\right)
  \;\Longrightarrow\;
  \left(\frac{4B}{A}\right)^r
  \leq
  p^{-1/2}
  %A\geq 4B\exp\left(\frac{1}{2B}\right)
$$
we obtain
\begin{eqnarray}
\nonumber
\frac{6}{\pi^2}\left(1-u^{-1}\right)H^2f(X)
&\leq&
4A^2(\log p)^2p^{1/2}
\left[
\frac{1}{A\log p}\left(\frac{4B}{A}\right)^rp^{1/2}+\frac{3B}{A}
\right]
\\
\nonumber
&\leq&
4A^2p^{1/2}(\log p)^2
\left(
\frac{1}{A\log p}+\frac{3B}{A}
\right)
\\
\label{E:explain2}
&=&
12ABp^{1/2}(\log p)^2
\left(
1+\frac{1}{3B\log p}
\right)
\,.
\end{eqnarray}

Plugging in our choices of $A$ and $B$ yields:~\footnote{At this point our choices of $A$ and $B$ are properly motivated --
the condition \mbox{$A\geq 4B\exp(1/(2B))$} was to ensure that
the quantity in the square brackets of (\ref{E:explain1}) remains bounded as $p\to\infty$, and we wanted to minimize $AB$ so that
the constant appearing in (\ref{E:explain2}) was as small as possible.}
\begin{equation}\label{E:subfinal}
\frac{6}{\pi^2}\left(1-u^{-1}\right)H^2f(X)
\leq
3e^2p^{1/2}(\log p)^2\left(1+\frac{4}{3\log p}\right)
\end{equation}
As $f(X)$ is increasing and positive for $X\geq 14$,
the result now follows upon solving
(\ref{E:subfinal}) for $H$.
\qed

\vspace{1ex}
\noindent\textbf{Proof of Theorem~\ref{T:1}.}
%We adopt all notation from Theorem~\ref{T:q2maintech}.
Suppose $p\geq 10^7$.
Let $n_0$ denote the smallest $n\in\Z^+$ such that $(n,u)=1$ and $\chi(n)\neq 1$.
Set $H:=n_0-1$ so that
$\chi(n)=1$ for all $n\in[1,H]$ with $(n,u)=1$.

First we show that 
%In order to apply Theorem~\ref{T:q2maintech}, it suffices to show that
$H\leq (2e^2\log p-2)^{1/2}p^{1/2}$.  By way of contradiction,
suppose $H>(2e^2\log p-2)^{1/2}p^{1/2}$.  In this case
we set $H_0=\lfloor (2e^2\log p-2)^{1/2}p^{1/2}\rfloor$,
and note that we still have $\chi(n)=1$ for all $n\in[1,H_0]$ with $(n,u)=1$ for this smaller value $H_0$.
We invoke Theorem~\ref{T:2} to conclude that
$H_0<Kg(p)p^{1/4}\log p$ where $Kg(p)\leq Kg(10^{7})< 12$. %11.05
Using again the fact that $p\geq 10^{7}$, we have
$$
H_0<12 p^{1/4}\log p<(2e^2\log p-2)^{1/2}p^{1/2}-1<H_0 %11.05
\,,
$$
which is
a contradiction.
This proves that
$H\leq (2e^2\log p-2)^{1/2}p^{1/2}$.
 
Having shown that $H$ satisfies the required condition,
we apply Theorem~\ref{T:2} to find
%Thus Theorem~\ref{T:q2maintech} applies to give
$H\leq Kg(p_0)\,p^{1/4}\log p$ when $p\geq p_0\geq 10^7$.
Therefore
$$
  n_0\leq Kg(p_0)\,p^{1/4}\log p+1
  \,,
$$
 for $p\geq p_0\geq 10^7$.  
Computation of the table of constants is routine;
for each value of $p_0$, we compute (being careful to round up)
the quantity
$$
K g(p_0)+\frac{1}{p_0^{1/4}\log p_0}
\,.
\;\;\text{\qed}
$$

%
%
%
%
%\begin{lemma}\label{L:q2weak}
%Let $\chi$ be a non-principal Dirichlet character modulo a prime $p$.
%If $q_1<q_2$ are the two smallest primes such that $\chi(q_1),\chi(q_2)\neq 1$,
%then
%$$
%  q_2<2\,p^{1/2}\log p
%  \,.
%$$
%\end{lemma}

\section{Proofs of the Corollaries}\label{S:corollaries}

\vspace{1ex}
\noindent\textbf{Proof of Corollary~\ref{C:1}.}
Apply Theorem~\ref{T:1} with $u=q_1$ and 
observe that the smallest $n\in\Z^+$ with $(n,q_1)=1$ and $\chi(n)\neq 1$ is equal to $q_2$.
\qed

\vspace{1ex}
The following is a lemma due to Hudson (see~\cite{hudson:1973})
that will allow us to prove Corollary~\ref{C:2}.
The proof is brief and so we include it for the sake of completeness.

\begin{lemma}[Hudson]\label{L:hudson}
Let $\chi$ be a non-principal Dirichlet character modulo a prime $p\geq 5$.
Suppose that $q_1<q_2$ are the two smallest prime non-residues of $\chi$,
%Suppose $q_1<q_2$ are the two smallest primes such that $\chi(q_1),\chi(q_2)\neq 1$,
and that $q_1\neq 2$ or $q_2\neq 3$.
Let $S$ denote the maximal number of consecutive integers for which $\chi$ takes the same value.
%that is, suppose that $\chi$ is not constant on any half-open interval of length~$S$.
Then $q_2\leq S q_1+1$.
\end{lemma}

\noindent\textbf{Proof.}
Let $t\in\Z^+$ be maximal such that $1+t q_1<q_2$.
(This is always possible unless $q_1=2$ and $q_2=3$.)
Then the $t+1$ integers
\begin{equation}\label{E:list}
  1,\,1+q_1,\,\dots,\,1+tq_1
\end{equation}
are residues with respect to $\chi$.
Let $x$ be denote the unique inverse of $q_1$ modulo $p$
in the interval $(0,p)$.
Multiplying (\ref{E:list}) by $x$ allows us to see that the $t+1$ consecutive integers
$$
  x,\,x+1,\,\dots,\,x+t
$$
take on the same character value; hence $t+1\leq S$.  By the maximality of $t$,
we conclude that $q_2\leq(t+1)q_1+1\leq S q_1+1$.
\qed

\vspace{1ex}
We note that the above Lemma can be improved if $\chi(-1)=1$ (see~\cite{hudson:1973})
but we will not require this.
The other result we we use in the proof of Corollary~\ref{C:2}
is the following, which is a special case of
Theorem~1.2 of~\cite{mcgown:consecutive}.
\begin{theorem}\label{T:consecutive}
If $\chi$ is any non-principal Dirichlet character to the prime modulus
%$p\geq 5\cdot 10^{18}$
$p\geq 10^{19}$
which is constant on $(N,N+H]$, then
$H<7.1\,p^{1/4}\log p$.
\end{theorem}

\vspace{1ex}
\noindent\textbf{Proof of Corollary~\ref{C:2}.}
If $q_1>e^2\log p$, then we apply Corollary~\ref{C:1} and we are done.
Hence we may assume that $q_1\leq e^2\log p$.
If $q_2=3$, then we are clearly done, so we may also assume $q_2\neq 3$.
In this case, we combine Lemma~\ref{L:hudson} and Theorem~\ref{T:consecutive}
to conclude that $q_2\leq (7.1\,p^{1/4}\log p)(e^2\log p)+1<53\,p^{1/4}(\log p)^2$.~\qed

\vspace{1ex}
In order to prove Corollary~\ref{C:3}, we will use the following result which gives a weak bound on $q_2$,
but requires no extra hypotheses on $q_1$.

\begin{lemma}\label{L:q2weaknew}
Let $\chi$ be a non-principal Dirichlet character modulo $m\geq 10^{15}$.
Suppose that $q_1<q_2$ are the two smallest prime non-residues of $\chi$.
%Suppose $q_1<q_2$ are the two smallest primes such that $\chi(q_1),\chi(q_2)\neq 1$.
Then
$$
  q_2<
  2\,m^{1/2}\log m
  \,.
$$    
%  \begin{cases}
%    3\,m^{1/2}\log m & \text{ if $m\geq 10^7$ }\\
%    2\,m^{1/2}\log m & \text{ if $m\geq 10^{15}$ }
%  \end{cases}
%  \,.
%$$
\end{lemma}

\noindent\textbf{Proof.}
Using the explicit version of the P\'olya--Vinogradov inequality proven in~\cite{bachman.rachakonda},
we find
\begin{eqnarray*}
  %\left|
  \vrule \
  %\sum_{\footnotesize{\begin{array}{cc}n<x\\[-1ex](n,q_1)=1\end{array}}}
  \sum_{\substack{n<x\\(n,q_1)=1}}^{\phantom{n}}
  \chi(n)
  \ \vrule
  %\right|
  &=&
  \left|
  \sum_{n<x}
  \chi(n)
  -
  \chi(q_1)\sum_{n<x/q_1}
  \chi(n)
  \right|
  \\
  &\leq&
  \left|
    \sum_{n<x}
  \chi(n)
  \right|
  +
  \left|
  \sum_{n<x/q_1}
  \chi(n)
  \right|
  \\[1ex]
  &\leq&
  2\left(\frac{1}{3\log 3}\,m^{1/2}\log m+6.5\,m^{1/2}\right)
%  \\[1ex]
%  &\leq&
%  m^{1/2}\log m
  \,.
\end{eqnarray*}
If $\chi(n)=1$ for all $n\leq x$ with $(n,q_1)=1$, then
$$
  %\left|
  \vrule \
  %\sum_{\footnotesize{\begin{array}{cc}n<x\\[-1ex](n,q_1)=1\end{array}}}
  \sum_{\substack{n<x\\(n,q_1)=1}}^{\phantom{n}}
  \chi(n)
  \ \vrule
  %\right|
  \geq
  (1-q_1^{-1})x-1
  \,.
$$
Thus for $1<x<q_2$, we have
$$
  (1-q_1^{-1})x-1
  \leq
    2\left(\frac{1}{3\log 3}\,m^{1/2}\log m+6.5\,m^{1/2}\right)
    \,.
$$
Using the fact that $q_1\geq 2$ and letting $x$ approach $q_2$ from the left, we obtain
$$
  q_2
  \leq
  4\left(\frac{1}{3\log 3}\,m^{1/2}\log m+6.5\,m^{1/2}\right)+2
  \,,
$$
and the result follows.
\qed

\vspace{1ex}
\noindent\textbf{Proof of Corollary~\ref{C:3}.}
If $q_1<e^2\log p$, we use Lemma~\ref{L:q2weaknew} to obtain
\mbox{$q_2<2\,p^{1/2}\log p$}
and hence $q_1 q_2<2e^2p^{1/2}(\log p)^2<15\,p^{1/2}(\log p)^2$.
If $q_1\geq e^2\log p$, then we apply Theorem~\ref{T:norton}
(using the fact that $\chi$ has odd order)
and Corollary~\ref{C:1} to find
$q_1 q_2\leq C'\,p^{1/2}(\log p)^2$ with
%$C'=(3.9)(6.1537)<24$.
%$C'=(3.9)(6.1495)<24$.\qed
$C'=(3.9)(6.1536)<24$. \qed

\bibliographystyle{model1-num-names}
\bibliography{myrefs}

%\bibliographystyle{plain}
%\bibliography{myrefs}

\end{document}